\documentclass[12pt,leqno]{article}
\usepackage{amsfonts}

\usepackage{fullpage}
\usepackage{amssymb}
\usepackage{amsmath}
\usepackage{amsthm}
\newtheorem{theorem}{Theorem}

\newtheorem{lemma}{Lemma}

\begin{document}
\title{\bf \ An identity on the $2m$-th power mean value of
the generalized Gauss sums\footnote{This work was supported by the
National Natural Science Foundation of China, Grant No.
11071121.}}
\author{Feng Liu\, and Quan-Hui Yang\thanks {Corresponding author.
E-mail addresses: fliu\_1986@126.com; yangquanhui01@163.com.}\\
\small    School of Mathematical Sciences, Nanjing Normal
University,\\
\small  Nanjing  210046,  P. R. China }
\date{}
\maketitle \baselineskip 18pt %\mathcal{}

{\bf Abstract.} In this paper, using combinatorial and analytic
methods, we prove an exact calculating formula on the $2m$-th
power mean value of the generalized quadratic Gauss sums for
$m\geq 2$. This solves a conjecture of He and Zhang [`On the
$2k$-th power mean value of the generalized quadratic Gauss sums',
Bull. Korean Math. Soc. 48 (2011), No.1, 9-15].

  \vskip 5mm

\section{Introduction}
Let $q\geq 2$ be an integer and $\chi$ be a Dirichlet character
modulo $q$. For any integer $n$, the classical quadratic Gauss
sums $G(n;q)$ and the generalized quadratic Gauss sums
$G(n,\chi;q)$ are defined respectively by
$$G(n;q)=\sum^{q}_{a=1}e\left(\frac{na^2}{q}\right),$$
and
$$G(n,\chi;q)=\sum^{q}_{a=1}\chi(a)e\left(\frac{na^2}{q}\right),$$
where $e(x)=e^{2\pi ix}$.

The study of $G(n,\chi;q)$ is important in number theory, since it
is a generation of $G(n,q)$. In \cite{weil}, Weil proved that if
$p\geq 3$ is a prime, then
$$|G(n,\chi;p)|\leq 2\sqrt{p}.$$ In fact, Cochrane and Zheng \cite{cochrane} generalized this result
to any integer. That is, for any integer $n$ with $(n,q)=1,$ we
have
$$|G(n,\chi;q)|\leq 2^\omega(q) \sqrt{q},$$ where $\omega(q)$ is the number of all distinct prime divisors of $q$.

Beside the upper bound of $G(n,\chi;q)$, the power mean value of
$|G(n,\chi;q)|$ had also been studied by some authors. W. Zhang
(see \cite{zhang02}) proved that if $p$ is an odd prime and $n$ is
an integer with $(n,p)=1$, then
$$\sum_{\chi \mod p}|G(n,\chi;p)|^4=\left\{ \begin{array}{ll}
 (p-1)[3p^2-6p-1+4\left(\frac{n}{p}\right)\sqrt{p}]~,& p\equiv 1 \mod~4~;\\
 ~~~~~~~(p-1)(3p^2-6p-1)~,& p\equiv 3 \mod~4~.\end{array}
 \right.$$
 and $$\sum_{\chi \mod p}|G(n,\chi;p)|^6=(p-1)(10p^3-25p^2-4p-1),~~\text{if}~~ p\equiv 3 \mod 4,$$
where $(\frac{n}{p})$ is the Legendre symbol. For $p\equiv 1 \mod
4$, it is still an open question to calculate the exact value of
$\sum_{\chi \mod p}|G(n,\chi;p)|^6$.

In 2005, W. Zhang and H. Liu \cite{zhangliu} proved that if $q\geq
 3$ is a square-full number, then
 for any integer $n,k$ with $(nk,q)=1,k\geq
 1$, we have
  $$\sum_{\chi \mod q}|G(n,\chi;q)|^4=q\cdot\phi^2(q)\underset{p\mid q}
  {\prod}(k,p-1)^2\cdot \prod_{\substack{p\mid q\\(k,p-1)=1}}\frac{\phi(p-1)}{p-1},$$ where
  $\phi(q)$ is the Euler funtion.

Recently, Y. He and W. Zhang \cite{hezhang} proved the following
result.

Let odd number $q > 1$ be a square-full number. Then for any
integer $n$ with $(n,q)=1$ and $k$=3 or 4, we have the identity
$$\sum_{\chi \mod q}{|G(n,\chi;q)|^{2k}}=4^{(k-1)\omega(q)}\cdot
q^{k-1} \cdot \phi^2{(q)}.$$

Besides, they conjectured the above identity also holds for $k\geq
5$.

In this paper, we prove this conjecture in the following.
\begin{theorem}\label{thm1} Let odd number $q > 1$ be a
square-full number, $m\geq 2$ be an integer. Then for any integer
$n$ with $(n,q)=1$, we have the identity
$$\sum_{\chi \mod q}{|G(n,\chi;q)|^{2m}}=4^{(m-1)\omega(q)}\cdot
q^{m-1} \cdot \phi^2{(q)}.$$
\end{theorem} \vskip 2mm

\section{Proofs}

Let $p\geq 3$ be a prime, and let $k,n,a$ be three integers with
$1\leq k \leq n$.
 Write $$T_p(n,k,a)=\underset{x_1+\cdots+x_n\equiv a \mod p}{\sum^{p-1}_{x_1=1}\sum^{p-1}_{x_2=1}\cdots\sum^{p-1}_{x_n=1}}
  \left(\frac{x_1x_2\cdots x_k}{p}\right).$$
In order to prove Theorem \ref{thm1}, we need some lemmas on the
value of $T_p(n,k,a)$.

\begin{lemma}\label{lemma3}(See {[4, Theorem 8.2]}.) Let $p\geq 3$ be a prime. Then for any integer $a$, we have
 $$\sum^{p-1}_{x=1} \left(\frac{x^2+ax}{p}\right)
 = \left\{ \begin{array}{ll} -1, & \text{if}~~p\nmid a;
 \\ p-1,& \text{if}~~p\mid a. \end{array} \right.$$
\end{lemma} \vskip 2mm

This is a basic lemma which we will use to calculate the value of
$T_p(n,k,a)$.
\begin{lemma}\label{lemma4} Let $p\geq 3$ be a prime. Then for any integer $n\geq 1$, we have
 $$T_p(n,n,0)=\left\{ \begin{array}{ll} ~~~~~~~~~~0~\mbox{,} & \text{if}~~2\nmid n\mbox{;}
 \\ p^{(n-2)/2}(p-1)\left(\frac{-1}{p}\right)^{n/2}~\mbox{,}& \text{if}~~2\mid n. \end{array} \right.$$
\end{lemma} \vskip 2mm

\begin{proof}
  By the definition of $T_p(n,k,a)$ and Lemma \ref{lemma3}, for $n\geq 3,$ we have
 \begin{eqnarray*}
 & &
 T_p(n,n,0)\\
 &=&
 \underset{x_1+\cdots+x_{n}\equiv 0 \mod p}{\sum^{p-1}_{x_1=1}\sum^{p-1}_{x_2=1}\cdots
 \sum^{p-1}_{x_{n-1}=1}\sum^{p-1}_{x_{n}=1}}\left(\frac{x_1x_2\cdots
 x_{n-1}x_{n}}{p}\right)\\
 &=&~~
 \sum^{p-1}_{x_1=1}\sum^{p-1}_{x_2=1}\cdots\sum^{p-1}_{x_{n-1}=1}\left(\frac{x_1x_2\cdots
 x_{n-1}(-x_1-\cdots-x_{n-1})}{p}\right)\\
 &=&
 \left(\frac{-1}{p}\right)\sum^{p-1}_{x_1=1}\sum^{p-1}_{x_2=1}\cdots\sum^{p-1}_{x_{n-2}=1}\left(\frac{x_1x_2\cdots
 x_{n-2}}{p}\right)\\
 & &
 ~~~~~~~~~~~~~~~~~\cdot \sum^{p-1}_{x_{n-1}=1}\left(\frac{x_{n-1}^2+(x_1+\cdots
 +x_{n-2})x_{n-1}}{p}\right)\\
 &=&
 \left(\frac{-1}{p}\right)\underset{x_1+\cdots+x_{n-2}\not\equiv0 \mod p}{\sum^{p-1}_{x_1=1}
 \sum^{p-1}_{x_2=1}\cdots\sum^{p-1}_{x_{n-2}=1}}\left(\frac{x_1x_2\cdots
 x_{n-2}}{p}\right)\cdot (-1)\\
 & &+\left(\frac{-1}{p}\right)\underset{x_1+\cdots+x_{n-2}\equiv0 \mod p}
 {\sum^{p-1}_{x_1=1}\sum^{p-1}_{x_2=1}\cdots\sum^{p-1}_{x_{n-2}=1}}\left(\frac{x_1x_2\cdots
 x_{n-2}}{p}\right)\cdot (p-1)\\
 &=&
 \left(\frac{-1}{p}\right)(-1)\ {\sum^{p-1}_{x_1=1}\sum^{p-1}_{x_2=1}\cdots\sum^{p-1}_{x_{n-2}=1}}\left(\frac{x_1x_2\cdots
 x_{n-2}}{p}\right)\\
 & &-\left(\frac{-1}{p}\right)(-1)\underset{x_1+\cdots+x_{n-2}\equiv0 \mod p}
 {\sum^{p-1}_{x_1=1}\sum^{p-1}_{x_2=1}\cdots\sum^{p-1}_{x_{n-2}=1}}\left(\frac{x_1x_2\cdots
 x_{n-2}}{p}\right) \\
 & &
 +\left(\frac{-1}{p}\right)(p-1)\underset{x_1+\cdots+x_{n-2}\equiv0 \mod p}{\sum^{p-1}_{x_1=1}\sum^{p-1}_{x_2=1}
 \cdots\sum^{p-1}_{x_{n-2}=1}}\left(\frac{x_1x_2\cdots
 x_{n-2}}{p}\right)\\
 &=&
 \left(\frac{-1}{p}\right)p\underset{x_1+\cdots+x_{n-2}\equiv 0 \mod p}{\sum^{p-1}_{x_1=1}\sum^{p-1}_{x_2=1}
 \cdots\sum^{p-1}_{x_{n-2}=1}}\left(\frac{x_1x_2\cdots
 x_{n-2}}{p}\right)\\
 &=&
 \left(\frac{-1}{p}\right)p\cdot T_p(n-2,n-2,0).
 \end{eqnarray*}
It is easy to calculate $T_p(1,1,0)$ and $T_p(2,2,0)$.
 \begin{eqnarray*}
 & &T_p(1,1,0)=\sum^{p-1}_{x=1}\left(\frac{x}{p} \right)=0~\mbox{,}\\
 & &T_p(2,2,0)=\underset{x_1+x_2\equiv 0 \mod p}{\sum^{p-1}_{x_1=1}\sum^{p-1}_{x_2=1}}\left(\frac{x_1x_2}{p}\right)
 =\sum^{p-1}_{x_1=1}\left(
 \frac{-x_1^2}{p}\right)=\left(\frac{-1}{p}\right)(p-1).
 \end{eqnarray*}
Therefore, we have
 \begin{eqnarray*}
 & &T_p(2k+1,2k+1,0)=\left(\left(\frac{-1}{p}\right)p\right)^kT_p(1,1,0)=0~\mbox{,}\\
 & &T_p(2k,2k,0)=\left(\left(\frac{-1}{p}\right)p\right)^{k-1}T_p(2,2,0)=\left(\frac{-1}{p}\right)^kp^{k-1}(p-1).
 \end{eqnarray*}
 This completes the proof of Lemma \ref{lemma4}.
 \end{proof}

\begin{lemma}\label{lemma5} Let $p\geq 3$ be a prime and $n\geq 1$ be an integer.
Then for any integer $a$ with $(a,p)=1$, we have
 $$T_p(n,n,a)=\left\{ \begin{array}{ll} \left( \frac{a}{p}\right)
p^{(n-1)/2} \left( \frac{-1}{p}\right)^{(n-1)/2}, &
\text{if}~~2\nmid n;
 \\-p^{(n-2)/2}\left(\frac{-1}{p}\right)^{n/2},& \text{if}~~2\mid n. \end{array} \right.$$
\end{lemma} \vskip 2mm

\begin{proof}
Since $(a,p)=1,$ we have
 \begin{eqnarray*}
 T_p(n,n,a)
 &=&
 \underset{x_1+\cdots+x_n\equiv a \mod p}{\sum^{p-1}_{x_1=1}\sum^{p-1}_{x_2=1}\cdots\sum^{p-1}_{x_n=1}}
 \left(\frac{x_1x_2\cdots x_n}{p}\right)\\
 &=&
 \underset{ax_1+\cdots+ax_n\equiv a \mod p}{\sum^{p-1}_{ax_1=1}\sum^{p-1}_{ax_2=1}\cdots\sum^{p-1}_{ax_n=1}}
 \left(\frac{ax_1ax_2\cdots ax_n}{p}\right)\\
 &=&
 \left(\frac{a}{p}\right)^n\underset{x_1+\cdots+x_n\equiv 1 \mod p}
 {\sum^{p-1}_{x_1=1}\sum^{p-1}_{x_2=1}\cdots\sum^{p-1}_{x_n=1}}
 \left(\frac{x_1x_2\cdots x_n}{p}\right)\\
 &=&
 \left(\frac{a}{p}\right)^nT_p(n,n,1).
 \end{eqnarray*}
 The calculation of $T_p(n,n,1)$ is very similar to that of
 $T_p(n,n,0)$ in Lemma \ref{lemma4}, so we directly give the result here.
 $$T_p(n,n,1)=\left\{ \begin{array}{ll}
 p^{(n-1)/2}\left( \frac{-1}{p}\right)^{(n-1)/2}\mbox{,} & \text{if}~~2\nmid
 n;
 \\-p^{(n-2)/2}\left(\frac{-1}{p}\right)^{n/2}\mbox{,}& \text{if}~~2\mid n. \end{array}
 \right.$$
Hence,
 $$T_p(n,n,a)=\left(\frac{a}{p}\right)^nT_p(n,n,1)=\left\{ \begin{array}{ll} \left( \frac{a}{p}\right)
 p^{(n-1)/2}\left( \frac{-1}{p}\right)^{(n-1)/2}, & \text{if}~~2\nmid n\mbox{;}
 \\-p^{(n-2)/2}\left(\frac{-1}{p}\right)^{n/2},& \text{if}~~2\mid n. \end{array}
 \right.$$
 This completes the proof of Lemma \ref{lemma5}.
 \end{proof}

\begin{lemma}\label{lemma6}
Let $p\geq 3$ be a prime, and let $k,n,a$ be three integers with
$1\leq k \leq n$. Then we have
 $$T_p(n,k,a)=\left\{ \begin{array}{ll}
 (-1)^{n-k}\left(\frac{a}{p}\right)p^{(k-1)/2}\left(\frac{-1}{p}\right)^{(k-1)/2}~\mbox{,}
 & \text{if}~~2\nmid k~\mbox{and}~p\nmid a;\\
 ~~~~~~~~~~0~\mbox{,}& \text{if}~~2\nmid k~\mbox{and}~p~|~a;\\
 (-1)^{n+1-k}\left(\frac{-1}{p}\right)^{k/2}p^{(k-2)/2}~\mbox{,}& \text{if}~~2\mid k~\mbox{and}~p\nmid a;\\
 (-1)^{n-k}(p-1)\left(\frac{-1}{p}\right)^{k/2}p^{(k-2)/2}~\mbox{,}& \text{if}~~2\mid k~\mbox{and}~p~|~a.
 \end{array} \right. \eqno{(1)}$$
\end{lemma} \vskip 2mm
\begin{proof}[Proof.] For $k\leq n-1$, we have
\begin{eqnarray*}
 & &T_p(n,k,a)\\
 &=&\underset{x_1+\cdots+x_n\equiv a \mod p}
 {\sum^{p-1}_{x_1=1}\sum^{p-1}_{x_2=1}\cdots\sum^{p-1}_{x_n=1}}\left(\frac{x_1x_2\cdots
 x_k}{p}\right)\\
 &=&\sum^{p-1}_{x_1=1}\sum^{p-1}_{x_2=1}\cdots\sum^{p-1}_{x_{n-1}=1}\left(\frac{x_1x_2\cdots x_k}{p}\right)
 -\underset{x_1+\cdots+x_{n-1}\equiv a \mod p}{\sum^{p-1}_{x_1=1}\sum^{p-1}_{x_2=1}\cdots\sum^{p-1}_{x_{n-1}=1}}
 \left(\frac{x_1x_2\cdots x_k}{p}\right)\\
 &=&-T_p(n-1,k,a).
 \end{eqnarray*}
 Then by induction on $n$ we have
 $$T_p(n,k,a)=(-1)^{n-k}T_p(k,k,a)$$
  for all $n\geq k$.
 By Lemma \ref{lemma4} and Lemma \ref{lemma5}, we obtain equation (1), which completes the
 proof of Lemma \ref{lemma6}.

 \end{proof}

\begin{lemma}\label{lemma1} Let $p\geq 3$ be a prime, $\alpha \geq 2$, $a$ and $n$ be three integers with
    $1\leq a \leq p^\alpha -1$ and $(n,p)=1$. If $p^{\alpha-1}\parallel {a^2-1}$, then we write
    $a=rp^{\alpha-1}+\varepsilon$, where $1\leq r \leq p-1$ and $\varepsilon=\pm 1$. Then we have
$$\underset{b=1}{\overset{p^{\alpha}}{{\sum}'}}
 e\left(\frac{nb^2(a^2-1)}{p^\alpha}\right)=\left\{ \begin{array}{lll}
 0~\mbox{,}& \text{if}~~p^{\alpha -1}\nmid a^2-1;\\
p^{\alpha-1}\big[\big(\frac{2\varepsilon
rn}{p}\big)G(1;p)-1\big]~\mbox{,}&
\text{if}~~p^{\alpha -1}\parallel a^2-1;\\
\phi(p^\alpha)~\mbox{,}&
 \text{if}~~p^\alpha ~|~a^2-1~.\end{array}\right.$$
\end{lemma} \vskip 2mm
\begin{proof}[Proof.] See the proof of Lemma 4 of \cite{hezhang}.

\end{proof}

\begin{lemma}\label{lemma2} Let $p\geq 3$ be a prime. Then for any two integers $n\geq 1$ and $a$, we have
 $$\underset{x_1+x_2+\cdots+x_n\equiv a \mod p}{\sum^{p-1}_{x_1=1}\sum^{p-1}_{x_2=1}\cdots\sum^{p-1}_{x_n=1}}
 1=\left\{ \begin{array}{ll}
\big((p-1)^n-(-1)^n\big)\big/p~\mbox{,}& \text{if}~~p\nmid a;\\
\big((p-1)^n+(p-1)(-1)^n\big)\big/p~\mbox{,}&
\text{if}~~p~|~a.\end{array} \right.$$
\end{lemma} \vskip 2mm

\begin{proof}[Proof.] \begin{eqnarray*}
 \underset{x_1+\cdots+x_n\equiv a \mod p}{\sum^{p-1}_{x_1=1}\sum^{p-1}_{x_2=1}\cdots\sum^{p-1}_{x_n=1}}1
 &=&\sum^{p-1}_{x_1=1}\sum^{p-1}_{x_2=1}\cdots\sum^{p-1}_{x_{n-1}=1}1
    -\underset{x_1+\cdots+x_{n-1}\equiv a \mod p}{\sum^{p-1}_{x_1=1}\sum^{p-1}_{x_2=1}\cdots\sum^{p-1}_{x_{n-1}=1}}1 \\
 &=&(p-1)^{n-1}-\underset{x_1+\cdots+x_{n-1}\equiv a \mod
 p}{\sum^{p-1}_{x_1=1}\sum^{p-1}_{x_2=1}\cdots\sum^{p-1}_{x_{n-1}=1}}1.
\end{eqnarray*}
Then by induction on $n$, we have
\begin{eqnarray*}
\underset{x_1+\cdots+x_n\equiv a \mod
p}{\sum^{p-1}_{x_1=1}\sum^{p-1}_{x_2=1}\cdots\sum^{p-1}_{x_n=1}}1
&=&\sum^{n-2}_{k=1}(-1)^{k+1}(p-1)^{n-k}+(-1)^{n-2}\underset{x_1+x_2\equiv
 a \mod p}{\sum_{x_1=1}^{p-1}\sum^{p-1}_{x_2=1}}1\\
 &=&\left\{ \begin{array}{ll}
\big((p-1)^n-(-1)^n\big)\big/p\mbox{,}& \text{if}~~p\nmid a;\\
\big((p-1)^n+(p-1)(-1)^n\big)\big/p\mbox{,}&
\text{if}~~p~|~a.\end{array}
 \right.
 \end{eqnarray*}
 This completes the proof of Lemma \ref{lemma2}.
 \end{proof}

 \begin{lemma}\label{lemma7}(See [1, Theorem 9.13].)
For any odd prime $p$, we have
$$G^2(1;p)=\left(\frac{-1}{p}\right)p~.$$
\end{lemma} \vskip 2mm

 \begin{lemma}\label{lemma8}(See [7, Lemma 6].)
Let $m,n\geq 2$ and $u$ be three integers with $(m,n)=1$ and
$(u,mn)=1$. Then for any character $\chi=\chi_1\chi_2$ with
$\chi_1$ mod $m$ and $\chi_2$ mod $n$, we have the identity
$$G(u,\chi;mn)=\chi_1(n)\chi_2(m)G(un,\chi_1;m)G(um,\chi_2;n).$$
\end{lemma} \vskip 2mm

 \begin{lemma}\label{lemma9}
Let $p\geq 3$ be a prime, $\alpha \geq 2, m\geq 2$ be two
integers~. Then for any integer $n$ with $(n,p)=1$, we have the
identity
$$\sum_{\chi \mod p^\alpha}{|G(n,\chi;p^\alpha)|^{2m}}=4^{(m-1)}\cdot
\phi^2{(p^\alpha)} \cdot p^{(m-1)\alpha}.$$
\end{lemma} \vskip 2mm

\begin{proof}[Proof.] By the definition of
$G(n,\chi;p^\alpha),$ we have
\begin{eqnarray*}
  |G(n,\chi;p^{\alpha})|^2
  &=&\underset{a=1}{\overset{p^{\alpha}}{{\sum}'}}\underset{b=1}{\overset{p^{\alpha}}{{\sum}'}}
  \chi(a)\overline{\chi(b)}e\left(\frac{n(a^2-b^2)}{p^\alpha}\right)\\
  &=&\underset{a=1}{\overset{p^{\alpha}}{{\sum}'}} \chi(a)\underset{b=1}{\overset{p^{\alpha}}{{\sum}'}}
 e\left(\frac{nb^2(a^2-1)}{p^\alpha}\right)~.
  \end{eqnarray*}
Hence, by this formula we have
\begin{eqnarray*}
 & &\sum_{\chi\mod p^\alpha}|G(n,\chi;p^{\alpha})|^{2m}\\
  &=&\sum_{\chi\mod p^\alpha}\underset{x_1=1}{\overset{p^{\alpha}}{{\sum}'}}
  \underset{x_2=1}{\overset{p^{\alpha}}{{\sum}'}}
  \cdots \underset{x_m=1}{\overset{p^{\alpha}}{{\sum}'}}
    \chi(x_1\cdots x_m)\prod_{i=1}^{m}\left(\underset{y_i=1}
    {\overset{p^{\alpha}}{{\sum}'}}e\left(\frac{ny_i^2(x_i^2-1)}{p^\alpha}\right)\right)\\
&=&\phi(p^\alpha)\underset{x_1x_2\cdots x_m\equiv 1 \mod
p^\alpha}{\underset{x_1=1}
 {\overset{p^{\alpha}}{{\sum}'}}\underset{x_2=1}{\overset{p^{\alpha}}{{\sum}'}}
     \cdots \underset{x_m=1}{\overset{p^{\alpha}}{{\sum}'}}}
    \prod_{i=1}^{m}\left(\underset{y_i=1}
    {\overset{p^{\alpha}}{{\sum}'}}e\left(\frac{ny_i^2(x_i^2-1)}{p^\alpha}\right)\right).
 \end{eqnarray*}

Then by Lemma \ref{lemma1} we have
$$\sum_{\chi\mod
p^\alpha}|G(n,\chi;p^{\alpha})|^{2m}=\phi(p^\alpha)\sum_{k=0}^{m}{
{m \choose k}A(m,k)}, \eqno{(2)}$$ where
$$A(m,k)=\underset{x_1x_2\cdots
x_m\equiv 1 \mod p^\alpha}
   {\underset{\substack{x_1=1\\p^{\alpha-1} \parallel x_1^2-1}}{\overset{p^{\alpha}}
   {{\sum}'}}\cdots\underset{\substack{x_k=1\\p^{\alpha-1} \parallel x_k^2-1}}{\overset{p^{\alpha}}{{\sum}'}}
   \underset{\substack{x_{k+1}=1\\p^\alpha \mid  x_{k+1}^2-1}}{\overset{p^{\alpha}}{{\sum}'}}
     \cdots \underset{\substack{x_m=1\\p^\alpha \mid  x_m^2-1}}{\overset{p^{\alpha}}{{\sum}'}}}
    \prod_{i=1}^{m}\left(\underset{y_i=1}
    {\overset{p^{\alpha}}{{\sum}'}}e\left(\frac{ny_i^2(x_i^2-1)}{p^\alpha}\right)\right).$$
Now, in order to prove Lemma \ref{lemma9}, we need to calculate
$A(m,k)$.

\begin{eqnarray*}
 & &A(m,k)\\
 &=&\underset{x_1x_2\cdots x_{m}\equiv 1 \mod p^\alpha}
   {\underset{\substack{x_1=1\\p^{\alpha-1} \parallel x_1^2-1}}{\overset{p^{\alpha}}
   {{\sum}'}}\cdots\underset{\substack{x_k=1\\p^{\alpha-1} \parallel x_k^2-1}}{\overset{p^{\alpha}}{{\sum}'}}
   \underset{\substack{x_{k+1}=1\\p^\alpha \mid  x_{k+1}^2-1}}{\overset{p^{\alpha}}{{\sum}'}}
     \cdots \underset{\substack{x_{m}=1\\p^\alpha \mid  x_{m}^2-1}}{\overset{p^{\alpha}}{{\sum}'}}}
    \prod_{i=1}^{m}\left(\underset{y_i=1}
    {\overset{p^{\alpha}}{{\sum}'}}e\left(\frac{ny_i^2(x_i^2-1)}{p^\alpha}\right)\right)\\
 &=&2\phi(p^\alpha)\underset{x_1x_2\cdots x_{m-1}\equiv 1 \mod p^\alpha}
   {\underset{\substack{x_1=1\\p^{\alpha-1} \parallel x_1^2-1}}{\overset{p^{\alpha}}
   {{\sum}'}}\cdots\underset{\substack{x_k=1\\p^{\alpha-1} \parallel x_k^2-1}}{\overset{p^{\alpha}}{{\sum}'}}
   \underset{\substack{x_{k+1}=1\\p^\alpha \mid  x_{k+1}^2-1}}{\overset{p^{\alpha}}{{\sum}'}}
     \cdots \underset{\substack{x_{m-1}=1\\p^\alpha \mid  x_{m-1}^2-1}}{\overset{p^{\alpha}}{{\sum}'}}}
    \prod_{i=1}^{m-1}\left(\underset{y_i=1}
    {\overset{p^{\alpha}}{{\sum}'}}e\left(\frac{ny_i^2(x_i^2-1)}{p^\alpha}\right)\right)\\
 &=&2\phi(p^\alpha)A(m-1,k).
 \end{eqnarray*}
Hence, by induction on $m$, we have
$$A(m,k)=2^{m-k}\phi^{m-k}(p^\alpha)A(k,k).\eqno{(3)}$$

Next, we shall calculate $A(k,k)$. By the definition, we have

$$A(k,k)=\underset{x_1x_2\cdots x_k\equiv 1 \mod p^\alpha}
   {\underset{\substack{x_1=1\\p^{\alpha-1} \parallel  x_1^2-1}}{\overset{p^{\alpha}}
   {{\sum}'}}\cdots \underset{\substack{x_k=1\\p^{\alpha-1} \parallel x_k^2-1}}{\overset{p^{\alpha}}{{\sum}'}}}
    \prod_{i=1}^{k}\left(\underset{y_i=1}
    {\overset{p^{\alpha}}{{\sum}'}}e\left(\frac{ny_i^2(x_i^2-1)}{p^\alpha}\right)\right).$$
Write $x_i=r_ip^{\alpha-1}+\varepsilon_i(1\leq r_i \leq
p-1,\varepsilon_i=\pm 1)$ for $i=1,2,\cdots,k$. Then by Lemma
\ref{lemma1}, we have
\begin {eqnarray*}
& &A(k,k)\\
&=&p^{k(\alpha-1)}\underset{\substack{\varepsilon_1 r_1+\varepsilon_2r_2+\cdots +\varepsilon_kr_k\equiv 0 \mod p\\
   \varepsilon_1\varepsilon_2\cdots\varepsilon_k=1}}
   {\underset{r_1=1}{\overset{p-1}
   {{\sum}}}\underset{r_2=1}{\overset{p-1}
   {{\sum}}}\cdots \underset{r_k=1}{\overset{p-1}{{\sum}}}}
    \prod_{i=1}^{k}\left(\left(\frac{2n\varepsilon_ir_i}{p}\right)G(1;p)-1\right)\\
          &=&p^{k(\alpha-1)}\underset{\substack{ r_1+r_2+\cdots +r_k\equiv 0 \mod p\\
   \varepsilon_1\varepsilon_2\cdots\varepsilon_k=1}}
   {\underset{r_1=1}{\overset{p-1}
   {{\sum}}}\underset{r_2=1}{\overset{p-1}
   {{\sum}}}\cdots \underset{r_k=1}{\overset{p-1}{{\sum}}}}
    \prod_{i=1}^{k}\left(\left(\frac{2nr_i}{p}\right)G(1;p)-1\right)\\
          &=&2^{k-1}p^{k(\alpha-1)}\underset{ r_1+r_2+\cdots +r_k\equiv 0 \mod p}
   {\underset{r_1=1}{\overset{p-1}
   {{\sum}}}\underset{r_2=1}{\overset{p-1}
   {{\sum}}}\cdots \underset{r_k=1}{\overset{p-1}{{\sum}}}}
    \prod_{i=1}^{k}\left(\left(\frac{2nr_i}{p}\right)G(1;p)-1\right)\\
   &=&2^{k-1}p^{k(\alpha-1)}\cdot\underset{ r_1+r_2+\cdots+r_k\equiv 0 \mod p}
   {\underset{r_1=1}{\overset{p-1}
   {{\sum}}}\underset{r_2=1}{\overset{p-1}
   {{\sum}}}\cdots \underset{r_k=1}{\overset{p-1}{{\sum}}}}
    \Bigg((-1)^k\\
    &&+\left.\sum^k_{j=1}(-1)^{k-j}{k \choose j}\left(\frac{2n}{p}\right)^jG^j(1;p)\left(\frac{r_1r_2\cdots
    r_j}{p}\right)\right).
\end {eqnarray*}
By Lemma \ref{lemma6} and Lemma \ref{lemma2}, the above equality
becomes
\begin {eqnarray*}
& &A(k,k)\\
&=&2^{k-1}p^{k(\alpha-1)}(-1)^k \left(\frac{1}{p}\left((p-1)^k+(p-1)(-1)^k\right)\right.\\
   &&+\sum^{\lfloor k/2 \rfloor}_{j=1}(-1)^{2j}{k \choose
   2j}\left(\frac{2n}{p}\right)^{2j}G^{2j}(1;p)
   (-1)^{k-2j}\left(\frac{-1}{p}\right)^jp^{j-1}(p-1)\bigg).
\end {eqnarray*}
By Lemma \ref{lemma7}, we have
\begin {eqnarray*}
& &A(k,k)\\
&=&2^{k-1}p^{k(\alpha-1)-1}\Big((-1)^k(p-1)^k+(p-1)
   +\sum^{\lfloor k/2 \rfloor}_{j=1}{k \choose 2j}p^{2j}(p-1)\Big)\\
   &=&2^{k-1}p^{k(\alpha-1)-1}\Big((-1)^k(p-1)^k+(p-1)\big((p+1)^k+(1-p)^k\big)\big/2\Big)\\
   &=&2^{k-2}p^{k(\alpha-1)-1}\Big((p+1)(1-p)^k+(p-1)(p+1)^k\Big).
\end {eqnarray*}

Hence, by (3) we have
$$A(m,k)=2^{m-2}p^{m(\alpha-1)-1}\Big((-1)^k(p+1)(p-1)^m+(p-1)^{m-k+1}(p+1)^k  \Big).$$
Finally, by (2) we have
\begin {eqnarray*}
& &\sum_{\chi\mod
p^\alpha}|G(n,\chi;p^{\alpha})|^{2m}\\
&=&\phi(p^\alpha)\sum_{k=0}^{m}{ {m \choose
k}2^{m-2}p^{m(\alpha-1)-1}\Big((-1)^k(p+1)(p-1)^m+(p-1)^{m-k+1}(p+1)^k
\Big)}\\
&=&\phi(p^\alpha)2^{m-2}p^{m(\alpha-1)-1}(p+1)(p-1)^m\sum^m_{k=0}{m \choose k}(-1)^k\\
 & &~+\phi(p^\alpha)2^{m-2}p^{m(\alpha-1)-1}\sum^m_{k=0}{m \choose
 k}(p-1)^{m-k+1}(p+1)^k\\
 &=&0+\phi(p^\alpha)2^{m-2}p^{m(\alpha-1)-1}(p-1)(2p)^m\\
 &=&4^{m-1}\phi^2(p^\alpha)p^{\alpha(m-1)}.
\end {eqnarray*}
This completes the proof of Lemma \ref{lemma9}.
\end{proof}

\begin{proof}[Proof of  Theorem 1.] Since $q$ is an odd square-full number, let $q=p_1^{\alpha_1}p_2^{\alpha_2}\cdots
  p_{\omega(q)}^{\alpha_{\omega(q)}}$, we have $\alpha_i\geq 2,i=1,\cdots
  \omega(q)$. For any integer $n$ with $(n,q)=1$, by Lemma \ref{lemma8} and Lemma \ref{lemma9}, we
  obtain
  \begin{eqnarray*}
  & &\sum\limits_{\chi \mod q}|G(n,\chi;q)|^{2m}\\
  &=&\prod^{\omega(q)}_{\substack{i=1\\p_i^{\alpha_i}\parallel q}}
    \sum_{\chi \mod p_i^{\alpha_i}}|G(nq/p_i^{\alpha_i},\chi;p_i^{\alpha_i})|^{2m}\\
  &=&\prod^{\omega(q)}_{\substack{i=1\\p_i^{\alpha_i}\parallel q}}
     \big(4^{m-1}p_i^{\alpha_i(m-1)}\phi^2(p_i^{\alpha_i}) \big)\\
  &=&4^{(m-1)\omega(q)}\cdot q^{m-1}\cdot\phi^2(q)~.
 \end{eqnarray*}
 This completes the proof of theorem \ref{thm1}.
\end{proof}

\end{document}